# Regularity of Area-Minimizing Surfaces in 3D Polytopes and of Invariant Hypersurfaces in $\mathbf{R}^n$

Frank Morgan
Department of Mathematics and Statistics
Williams College
Williamstown, MA 01267
Frank.Morgan@williams.edu

## Abstract

In (the surface of) a convex polytope $P^3$ in $\mathbf{R}^4$, an area-minimizing surface avoids the vertices of P and crosses the edges orthogonally.

In a smooth Riemannian manifold M with a group of isometries G, an area-minimizing G-invariant oriented hypersurface is smooth (except for a very small singular set in high dimensions). Already in 3D, area-minimizing G-invariant *unoriented* surfaces can have certain singularities, such as three orthogonal sheets meeting at a point. We also treat other categories of surfaces such as rectifiable currents modulo and soap films.

# 1. Introduction

Standard geometric regularity theory deals with minimal objects in smooth manifolds. Yet applications from material science to cosmology require singular ambients or imposed symmetries. This paper studies such contexts.

Section 3 considers singular ambients modeled on 3D polytopes in $\mathbf{R}^4$, for example the surface of the hypercube or 3D manifolds with similar singularities, and proves that an area-minimizing oriented surface avoids the vertices and crosses the edges orthogonally. As a corollary we recover the regularity theorem of B. White [HPRR, Thm. 3] for isoperimetric oriented surfaces with imposed symmetry in a smooth manifold, in the case that the symmetries are locally orientation preserving, as follows by passing to the singular quotient. Recent mathematical study of such objects ([R1], see [R2, Sect. 1.6]) has application to crystallography (e.g. [LoM]; see [R1] and the references therein); natural structures seem to minimize area with imposed periodicity or symmetry. Our Theorem 4.1 provides extensions to higher dimensions and to boundary regularity.

For area-minimizing invariant *unoriented* surfaces, singularities can occur, such as three sheets meeting orthogonally at a point with cubical symmetry. Our Theorem 5.3 classifies the singularities in 3D; see Figure 5.1.

If isometries which do not preserve the orientation of the ambient manifold are admitted, then surfaces with orientation provided by oriented tangent planes rather than by choice of unit normal tend to admit more symmetries. We consider such surfaces (5.4—5.5) as well as soap films (5.11) and rectifiable currents modulo  (5.6—5.10). New singularities include cones over regular polyhedra for example.

The presence of interior singular curves leads to Cantor sets of boundary singularities (Rmk. 5.3.2, [M2, Fig. 13.9.3]). Even without imposed symmetries, the classification of area-minimizing tangent cones at the boundary remains open for soap films and rectifiable currents modulo  .

**1.1. The proofs.** The proofs begin with analysis of the possible area-minimizing tangent cones. The presence of ambient singularities or imposed symmetry eliminates some possible tangent cones but adds many new ones. Some can be ruled out by instability arguments [M1].

In showing that the surface is locally diffeomorphic to its tangent cone, imposed symmetry can often be used to simplify the proof. For example, when k sheets meet with imposed $\mathbf{Z}_k$ symmetry along a singular curve, the curve is automatically a smooth geodesic as the set of fixed points of the symmetry group.

**1.2. Acknowledgments.** I would like to thank Antonio Ros and the Department of Geometry and Topology at the University of Granada for inspiration and hospitality during my visit in May and June, 2003. This research is partially supported by a National Science Foundation grant.

# 2. Preliminaries

**2.1. Polytopes and singular ambient manifolds, rectifiable currents.** For this paper, a *polytope* $P^n$ in $\mathbf{R}^{n+1}$ is the boundary of a compact, positive-volume intersection of finitely many halfspaces. A polytope $P^n$ is a Riemannian manifold with codimension-two singularities, which do not affect the computation of n-volume or (n–1)-area. For example, the tetrahedron has singular vertices, but the edges are not intrinsic singularities.

As a "compact Lipschitz neighborhood retract," $P^n$ inherits from $\mathbf{R}^{n+1}$ the theory of *rectifiable currents,* the possibly singular oriented surfaces of geometric measure theory (see [Fed1], or [M2], especially sects. 4.3, 5.5, 9.1). The theory applies to any manifold $M^n$ in $\mathbf{R}^N$ locally (smoothly) diffeomorphic to polytopes. If M is compact, the Compactness Theorem guarantees the existence of an area-minimizing rectifiable current with given boundary or homology. One may impose the volume constraint that for any competitor + R, the algebraic volume of R is 0. Or one may impose invariance under some group G of isometries of M.

**2.2. Unoriented surfaces** [M2, Chapt. 11]. In geometric measure theory, unoriented surfaces are defined as rectifiable currents modulo two, and more generally one considers rectifiable currents modulo any integer 2. The previous section 2.1 (except for volume constraints) applies as well to these spaces of surfaces.

**2.3. Tangent cones.** Let be an area-minimizing (or merely stationary or bounded-mean-curvature) surface in $\mathbf{R}^n$. Then at any interior point or point of smooth prescribed boundary, has a well-defined density and at least one tangent cone C (weak limit of homothetic expansions) in the tangent space ([A1, 6.5], [A2, Sect. 3]; see [M2, 9.7]). The tangent cone C has the same density as .

These results apply to minimizers with prescribed volume, because they have (weakly) bounded mean curvature [M2, 13.5]. If is area minimizing (with or without prescribed volume), C is minimizing without prescribed volume (since volume becomes relatively negligible under large homothetic expansions).

These results and proofs extend immediately to polytopes $P^n$, and then as we will now explain to manifolds M in $\mathbf{R}^N$ locally diffeomorphic to polytopes. Suppose that $0 \in M$, and let F be a local diffeomorphism of $\mathbf{R}^N$ carrying M to a polytope P. We may assume that F is the identity to first order. The Riemannian metric that F induces on M well approximates the standard metric on M. As explained in [M4, Prop. 3.2], F( ) has a certain "monotonicity" property which implies estimates on the area of inside small balls about 0 and hence the existence of density and a tangent cone.

**2.4. Allard's regularity theorems.** Allard's regularity theorems apply in the very general context of varifolds in $\mathbf{R}^N$, including rectifiable currents and flat chains modulo . They apply to surfaces which are minimizing or merely stationary or of (weakly) bounded mean curvature; for hypersurfaces this means that for smooth variations, area A and volume V satisfy

$$|dA/dt| \leq C\,|dV/dt|.$$

The theorems apply at an interior point where the density is near the local minimum, which for our cases of interest is 1, and at a boundary point where the density is near 1/2.

Allard's regularity theorems ([A1, Sect. 8], [A2, Sect. 4]) say that if in some ball away from the boundary the density ratio is near 1, or if in some ball about a smooth ($C^\infty$) piece of boundary the density ratio is near 1/2, then in a slightly shrunken ball the surface is the graph of a $C^1$ function f (and Df has small Hölder constant).

It follows that if away from boundary such a surface is weakly close to a plane, then it is $C^1$ close.

The theorems extend from Euclidean space to a smooth Riemannian manifold by isometrically smoothly embedding the manifold in Euclidean space (see [A1, Rmk. 4.4]).

Once a surface of zero or constant mean curvature is known to be $C^1$, higher smoothness ($C^\infty$) follows by Schauder theory ([My, 1.10.4(i), 1.11.3]; see [M4, Prop. 3.5]), assuming that the ambient and prescribed boundary are smooth.

**2.5. Maximum principle** ([Ho1], [Ho2]; see [Ser]). Hopf's maximum principle implies that two smooth hypersurfaces of the same constant mean curvature, one locally above (or equal to) the other, which coincide at a point coincide locally. (If they are area minimizing and one is smooth, the other is automatically smooth at a point of contact, where the tangent cone, contained in a halfplane, must be a plane.) Similarly Hopf's boundary maximum principle implies that two such hypersurfaces tangent at a smooth boundary point coincide locally.

# 3. Regularity in 3D Polytopes

Chapter 3 proves regularity for area-minimizing surfaces in 3D polytopes, or in manifolds with similar singularities.

**3.1. Theorem.** *In a polytope $P^3$ in $\mathbf{R}^4$, let $\Sigma$ be a surface which minimizes area with prescribed boundary or volume, among oriented or unoriented surfaces. Then on the interior $\Sigma$ is a smooth minimal or constant-mean-curvature surface which avoids the vertices of P and crosses the edges differentiably and orthogonally.*

*Remark.* For example, for any polytope $P^3$, for small prescribed volume, geodesic balls about some vertex are area minimizing [M3, Thm. 3.3].

*Proof.* Since unoriented surfaces are locally orientable on the interior, we may assume that    is an oriented surface (rectifiable current). Using a standard decomposition argument ([Fed1, 4.5.17], [M2, Lemma Sect. 10.1]), we may assume that    locally is the boundary of a nested region (of multiplicity one). Consider a tangent cone C to    at a point p. C is area minimizing. Its link    in the unit geodesic sphere S, which is stationary and has minimizing tangent cones, must consist of geodesics which avoid the singularities [CFG, Cor. 2.5].

If p is a regular point of P, then    is well known to be regular ([Fl], [M4]).

If p is a vertex of P, then S is locally isometric to a sphere of curvature $K_0 > 1$, except for finitely many conical singularities    (of positive curvature type) corresponding to edges of P emanating from p. Therefore S is unstable by [M1, Thm. 3.1]. Actually, [M1, Thm. 3.1] considers the case without singularities   , but the proof is local and applies to our case as well.

If p is on the interior of an edge of P, then S is locally isometric to a unit sphere, except for two antipodal conical singularities ("the poles").

The only possible link is the equator, and C is a plane orthogonal to the edge. By Allard's regularity theorem 2.4, near p, S is a $C^1$ surface $C^1$ close to C. It follows that at p, is a $C^1$ surface normal to the edge.

If due to our initial decomposition there are two sheets at a point, they must coincide by a maximum principle 2.5, unless they intersect only at a point p on an edge. We may assume that the edge is vertical and consider a small, narrow, vertical cylinder about p. Let $_0$ above $_1$ be the two nearly horizontal sheets inside the cylinder, with say upward unit normals. Translate $_1$ upward a small distance $> 0$ to obtain $_2$, with $_2$ still below $_0$. (Any enclosed volume has been increased by A, where A is the cross-sectional area of the narrow cylinder.) By interior regularity, $_0$ + $_2$ cannot be area minimizing in the cylinder. Take a smaller surface in the cylinder with the same boundary and volume, and decompose it into pieces $_0'$ and $_2'$ with the same boundaries as $_0$ and $_2$. Translate $_0'$ downward distance to obtain a surface $_1'$, restoring enclosed volume to the original value. Now $_0'$ + $_1'$ has the same boundary and volume as the original $_0$ + $_1$, but less area, a contradiction.

**3.2. Singular manifolds.** Theorem 3.1 and proof hold for a homologically area-minimizing surface, perhaps with prescribed boundary or volume, in a 3-manifold in $\mathbf{R}^N$ locally diffeomorphic to polytopes, with (convex) polytopal tangent cones, except that I do not know how to eliminate distinct sheets touching at a point on a singular curve in M (the missing ingredient is translation). This issue does not arise for boundaries of regions of multiplicity 1.

Manifolds with milder singularities, such as the surface of a cylindrical can, with a Lipschitz Riemannian metric, satisfy the standard regularity theory [M4].

**3.3. Remark.** The Regularity Theorem 3.1 fails in "nonconvex polytopes," even for polyhedra $P^2$ in $\mathbf{R}^3$: a length-minimizing curve can have many arcs meeting at a vertex that has vertex angle $>> 2$ [CFG, Fig.2]. More generally, Morgan [M1] considers $\mathbf{R}^n$ with metric

$$ds^2 = dr^2 + {}^2 r^2 d\ {}^2,$$

where d is the standard metric on the unit sphere, with an isolated singularity at the origin. For $\mathbf{R}^4$ with $4/3$, the cone over $\mathbf{S}^1 \times \mathbf{S}^1$ is area minimizing [M1, Prop. 3.3].

More generally, in the context of soap films or minimizing separators of regions, the plane, the Y, and the tetrahedral cone in $\mathbf{R}^3$ lift via the Hopf fibration to minimizing cones in $\mathbf{R}^4$ with $\leq 4/3$ (see [B2, Sect. 9]), minimizing because the associated "paired calibrations" on $\mathbf{R}^3$ ([LM, Sect. 1.1], [B1, Thm. 3.1]) lift to $\mathbf{R}^4$. (The plane lifts to the previously mentioned example of the cone over $\mathbf{S}^1 \times \mathbf{S}^1$.)

# 4. Invariant Oriented Hypersurfaces

Theorem 4.1 provides interior and boundary regularity for area-minimizing invariant oriented hypersurfaces, generalizing the 3D result of B. White [HPRR, Thm. 3].

**4.1. Theorem (Regularity for invariant oriented hypersurfaces).** *In a smooth ($C^\infty$) Riemannian manifold M, among hypersurfaces oriented by unit normal invariant under a group G of isometries of M, perhaps with given boundary or homology or volume, suppose that $\Sigma$ minimizes area. Then on the interior $\Sigma$ is a smooth constant-mean-curvature hypersurface except for a singular set of codimension at least 7 in $\Sigma$.*

*Along a smooth piece of boundary, $\Sigma$ locally decomposes into smooth submanifolds with boundary (pairwise disjoint or coincident).*

*4.1.1. Remarks.* For $G = I$, interior regularity is due to Federer ([Fed2, Thm. 1], see [M2, Chapt. 10]), and boundary regularity is due to Hardt and Simon [HaS] and White [W1]. Along a boundary $\mu B$ (with integer multiplicity $\mu \geq 1$), there may be components with boundary $-B$, but each is paired smoothly with a component with boundary $+B$. With nontrivial imposed symmetry, all components have boundary $+B$.

If the isometries locally preserve orientation, it makes no difference whether the hypersurfaces are oriented by unit normal or by oriented tangent planes, but we wish also to allow reflectional symmetry of a sphere in $\mathbf{R}^3$ for example. In particular, our results apply to isoperimetric surfaces invariant under isometries which leave the regions they bound invariant.

**4.2. Lemma.** *Let $\Sigma$ be minimizing (as in Theorem 4.1). Then $\Sigma$ has constant mean curvature (weakly).*

*Proof.* For the mean curvature to be weakly equal to a constant $H_0$ means that for any smooth vectorfield, $dA/dV = H_0$. By standard variational arguments this holds for any invariant vectorfield. But on an invariant surface, the effect of any vectorfield is the same as the effect of its average over G.

**4.3. Lemma.** *In any dimension, there is a constant $\Upsilon > 1$, such that if $\Sigma$ is a minimizing hypersurface (as in Theorem 4.1), then at any interior point of density less than $\Upsilon$, $\Sigma$ is a smooth ($C^\infty$) constant-mean-curvature submanifold. Also at any smooth-boundary point of density less than $\Upsilon/2$, $\Sigma$ is a smooth submanifold with boundary.*

*Remark.* Lemmas 4.2 and 4.3 apply as well to unoriented surfaces and more general varifolds. Without a volume constraint, "constant-mean-curvature" becomes "minimal."

*Proof of Lemma 4.3.* By Lemma 4.2, $\Sigma$ has constant mean curvature (weakly). By Allard's regularity theorems 2.4, $\Sigma$ is smooth.

*Proof of Theorem 4.1.* For interior regularity, given Lemma 4.3, the proof (by induction) is the same as for the case $G = I$ [Fed2, Thm. 1]. In particular, for a G-invariant hypersurface, the instabilities of Simons (the $f$ of [Sim, Lemma 6.1.7] or the $\varphi$ of [Fed1, 5.4.13]) are all G invariant, depending only on the curvature of the link of the tangent cone. Also, the standard local decomposition into boundaries of nested regions of multiplicity one is G invariant.

For the base case of a curve $\Sigma$ in a two-dimensional manifold, consider a tangent cone C at a point with nontrivial isotropy subgroup $G_1$. The tangent cone C consists of oriented rays and is minimizing without volume constraint. $G_1$ must be cyclic, dihedral, or generated by a single reflection. If $G_1$ is cyclic, the tangent space mod $G_1$ is a cone with angle less than $2\pi$ and C is not minimizing. If $G_1$ is a dihedral group the tangent space mod $G_1$ is a sector less than $\pi$ with free boundary, and the cone is not minimizing. (A ray of C cannot lie in the boundary of the sector, because $G_1$ reverses the orientation of its lift.) If $G_1$ is generated by a single reflection, the tangent space mod $G_1$ is a sector of angle $\pi$ with free boundary, C mod $G_1$ must be a single ray in the middle, C is a line, and $\Sigma$ is regular by Lemma 4.3.

At a point p in a smooth piece of boundary, for a new type of singularity, there must be imposed symmetry preserving the boundary; the only possibility is some nontrivial $\mathbf{Z}_k$. The boundary B must be the totally geodesic set of fixed points. Locally pass to the quotient orbifold M', where the surface $\Sigma'$ is area minimizing without imposed symmetry, with boundary B' with integer multiplicity $\mu \geq 1$. By a theorem of White [W1, p. 294], $\Sigma'$ has a local mass decomposition into $\mu$ area-minimizing surfaces $\Sigma_i$ with boundary B'. Fix one $\Sigma_i$. The argument of Hardt and Simon [HS, Step II of the proof of 11.1, pp. 477-479] shows that the tangent cone Q to $\Sigma_i$ must lie in a plane. In our orbifold, the only possibility for Q is a halfplane, which lifts locally to k halfplanes $H_j$ in the tangent space of M, symmetric under $\mathbf{Z}_k$.

> To summarize for the reader's convenience the argument of Hardt and Simon, as simplified in our case, consider the component $Q_0$ of the regular points of Q with boundary B and the angular polar coordinate $\theta$ from the 2-plane orthogonal to the plane tangent to B. Because $Q_0$ is an oriented surface with boundary B, $\theta$ has a single-valued branch on $Q_0$, which on the intersection of $Q_0$ with the unit sphere **S** attains a minimum value $\theta_0$. Now $Q_0$ must coincide locally with the halfplane $\theta = \theta_0$ by the maximum principle 2.5, applied either on the interior of $Q_0 \cap$ **S**, or at the boundary $B \cap$ **S**, where $Q_0 \cap$ **S** is regular by induction.
> To show that $Q = Q_0$, assume that $R = Q - Q_0 \neq 0$. Let $Q_1$ be the halfplane with boundary B antipodal to $Q_0$, let $W_+$, $W_-$ be the two components of the complement of $Q_0 \cup Q_1$, and let $R_\pm = R \cap W_\pm$. By interior regularity, $(\partial R_\pm) \cap Q_0 = 0$. Now rotate $Q_0$ until it first touches R, to eventually conclude by the maximum principle 2.5 that R consists of halfplanes with boundary B. Since R has no boundary and $R + Q_0$ is area minimizing, the only possibility in $\mathbf{R}^n$ would be that R is a plane containing $Q_0$; in the orbifold, there are no possibilities.

The surface $\Sigma_i$ in M′ lifts to a surface $\tilde\Sigma$ [%Tau] in M, with boundary kB. Again by the theorem of White [W1, p. 294], $\tilde\Sigma$ has a mass decomposition into k surfaces $\tilde\Sigma_j$ with boundary B, with tangent cone contained in $\cup H_j$. Hence each $\tilde\Sigma_j$ must have a halfplane tangent cone and density 1/2. By Allard's regularity theorem 2.4, $\tilde\Sigma_j$ is a smooth manifold with boundary. There are such contributions for each $\Sigma_i$. If two are tangent at p, we may assume by reassembly that one lies above the other, and hence by Hopf's boundary maximum principle 2.5 that they coincide.

# 5. Invariant Unoriented Surfaces in $\mathbf{R}^3$

Theorem 4.1 proved regularity for an area-minimizing surface oriented by unit normal, the natural setting for the isoperimetric problem. Within different classes of orientability, singularities can occur already in 2D. This chapter treats regularity in 3D for surfaces with oriented tangent planes (Theorem 5.5), for unoriented surfaces or rectifiable currents modulo two (Theorem 5.3), for rectifiable currents modulo 3 (5.6—5.10), and for soap films (5.11).

We begin with two propositions on length-minimizing curves in 2D manifolds. Proposition 5.1 notes how the standard regularity generalizes to imposed symmetry.

**5.1. Proposition (General invariant curve regularity).** *In a smooth ($C^\infty$) 2D Riemannian manifold M, a curve Γ with prescribed finite boundary that is length minimizing with imposed symmetry, or more generally that is stationary, consists of geodesics meeting at isolated points.*

*Proof.* A minimizer with imposed symmetry is indeed stationary (without imposed symmetry) by the simple argument of Lemma 4.2. The proposition follows from the more general structure theorem of Allard and Almgren [AA, Sect. 3]. We sketch a less technical proof.

Away from the finite set of boundary points and fixed points of the symmetry group, is an immersed geodesic.

At a general point, a tangent cone consists of rays from the origin. In a small annulus about the point, consists of nearly radial geodesics. The result follows.

**5.2. Proposition (Regularity for invariant unoriented curves).** *In a smooth ($C^\infty$) 2D Riemannian manifold M, among unoriented curves (flat chains modulo two) invariant under a group G of isometries of M, with given boundary or homology, suppose that Γ minimizes length. Then on the interior Γ consists of geodesics, k of which may cross at a point with imposed reflectional (k=2), $Z_{2k}$ (cyclic), or $D_k$ (dihedral) symmetry. At an isolated boundary point, an odd number k of geodesics may meet with imposed $Z_k$ symmetry, or three may meet in a T with reflectional symmetry. All such singularities can occur.*

*Proof.* At an interior singular point with (nontrivial) isotropy subgroup $G_1$, the tangent cone C consists of 2k  4 rays. $G_1$ must be cyclic, dihedral, or generated by a single reflection. If $G_1$ is cyclic, the tangent space mod $G_1$ is a cone with angle less than $2\pi$, C mod $G_1$ must be a single ray, and by Proposition 5.1 consists of k geodesics crossing with imposed $\mathbf{Z}_{2k}$ symmetry. If $G_1$ is a dihedral group [or generated by a single reflection], the tangent space mod $G_1$ is a sector with angle less than [or equal to] $\pi$ and free boundary rays of multiplicity 1/2. Except for the case of reflection, C mod $G_1$ must consist of one or both boundary rays, and consists of k geodesics crossing with imposed $D_k$ (or $D_{2k} \supset D_k$) symmetry. For reflection, C mod $G_1$ could in addition or instead contain a ray orthogonal to the boundary rays, and consists of a single geodesic or two geodesics crossing with reflectional symmetry.

The case of a boundary point is similar but more limited because k must be odd.

Consideration of the quotient orbifold shows every such singular cone in $\mathbf{R}^2$ minimizing.

Our main Theorem 5.3 characterizes singularities in invariant unoriented area-minimizing surfaces in 3D, as pictured in Figure 5.1. Note that if only orientation-preserving isometries are allowed, only singularities (a), (b) with k=2, and (f) can occur.

**5.3. Theorem (Regularity for invariant unoriented surfaces).** *In a smooth ($C^\infty$) 3D Riemannian manifold M, among unoriented surfaces (flat chains modulo two) invariant under a group G of isometries of M, with given boundary or homology, suppose that $\Sigma$ minimizes area.*

*Then on the interior $\Sigma$ is a smoothly immersed minimal surface with only the following singularities of Figure 5.1:*

*(a) $k \geq 2$ "vertical" sheets cross at equal angles along a smooth curve and G contains $Z_{2k}$ (cyclic) or $D_k$ (dihedral) or reflection (k=2 only);*

*(b) case (a) plus a smooth horizontal sheet and an additional reflective symmetry or just $D_k$ or if k=2 all orientation-preserving symmetries of the cube (the octahedral group);*

*(c) six sheets cross as in the cone over the 1-skeleton of Figure 5.1(c), and G contains the $S_4$ of symmetries of the tetrahedron;*

*(d) nine sheets cross as in the cone over the 1-skeleton of Figure 5.1(d), and G contains the symmetries of the cube, generated by reflections across axis and diagonal planes;*

*(e) fifteen sheets cross as in the cone over the 1-skeleton of Figure 5.1(e), and G contains the symmetries of the dodecahedron.*

*Along a smooth piece of boundary, $\Sigma$ consists of*

*(f) an odd number $k \geq 1$ of halfsheets meeting smoothly at equal angles and G contains $Z_k$,*

*(g) a halfsheet of reflective symmetry with boundary and an orthogonal sheet tangent to the curve, or*

*(h) case (f) plus a totally geodesic sheet of reflection orthogonal to the boundary curve, and G also contains the reflection.*

*All such singularities can occur.*

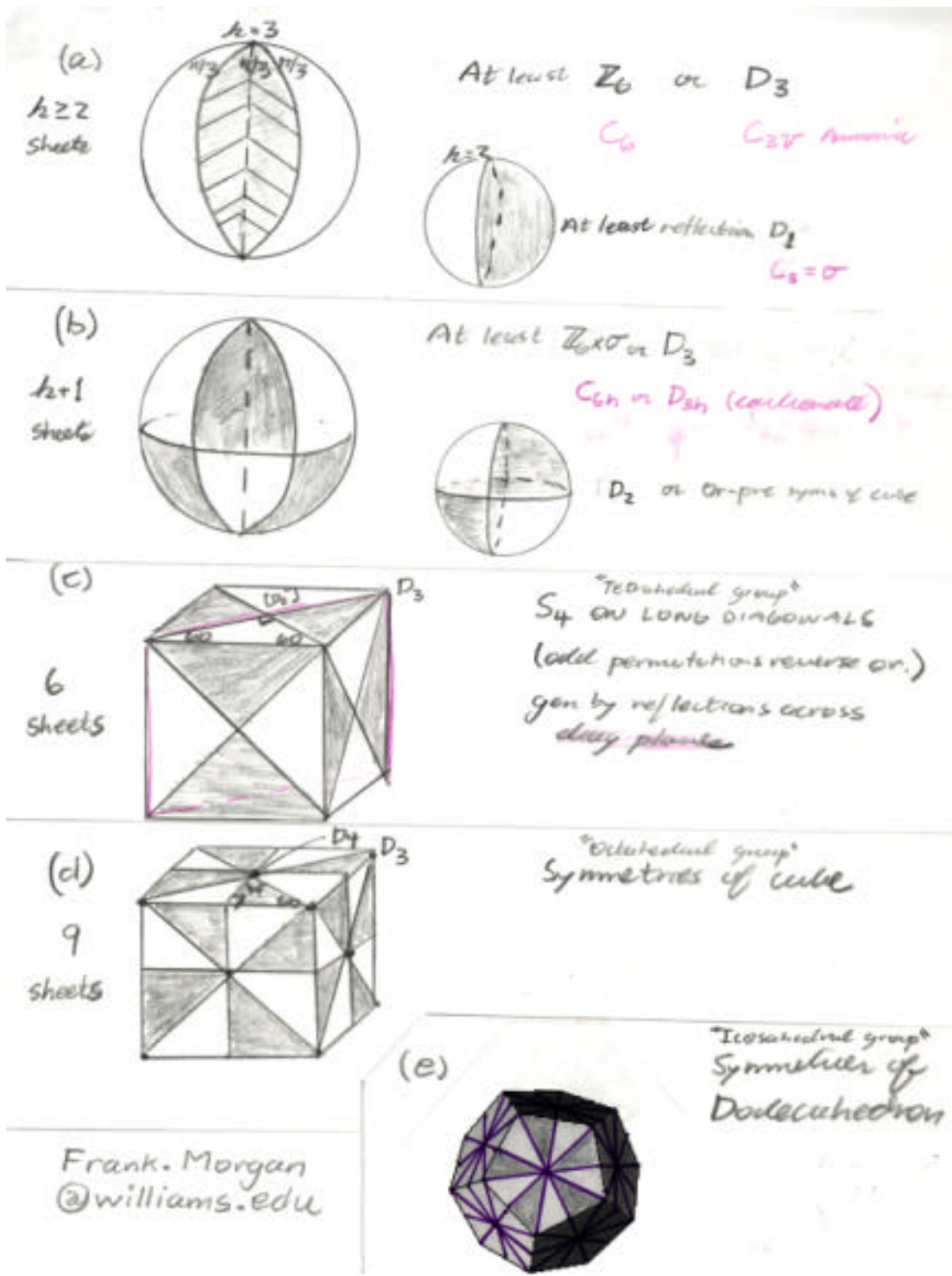

**Figure 5.1**
The five types of interior singularities for unoriented area-minimizing surfaces in 3D with imposed symmetries.

*5.3.1. Cubic torus.* A. Ros mentioned to me an interesting example of three orthogonal sheets (case (b), k=2). In a cubic 3D torus the minimizer with imposed axis reflectional symmetry in homology class (1,1,1) is three orthogonal flat 2D tori. (In the quotient orbifold there are no other possibilities.). The best oriented surface (and probably the best unoriented surface) without imposed symmetry is a totally geodesic torus represented by the hexagon through the middle of the cube, perpendicular to a long diagonal. (For an oriented minimizer , stable implies totally geodesic, hence a torus.)

*5.3.2. Remark.* Two sheets crossing with reflective symmetry (case (a)) admit a smooth symmetric boundary curve inside the surface of reflection which intersects the other sheet in a Cantor set (case (g)). A similar phenomenon occurs with surfaces with oriented tangent planes (5.5(d), Figure 5.3), with surfaces modulo three (5.7(a,b)), and with soap films (5.11).

*5.3.3. Isoperimetric problem.* Fischer and Koch [FiK] introduced (smoothly embedded) minimal "balance" surfaces which divide the ambient manifold M into two regions interchanged by some and possibly preserved by other imposed symmetries. Theorem 5.3 applies to area-minimizing balance surfaces with singularities. All cases (a)—(h) so occur, as indicated by the shadings in Figure 5.1. The proof goes through unchanged, since such surfaces are stationary and have minimizing tangent cones.

*Proof of Theorem 5.3. Part I: at an interior point p, a tangent cone C is a plane or of type (a)—(e).* The link in the unit sphere is stationary with minimizing tangent cones. By Propositions 5.1 and 5.2, the link consists of geodesics meeting in even numbers at equal angles. If there are no such junctions, C is a plane. If geodesics meet only in fours, the link consists of orthogonal geodesics, and the cone is of type (a) or (b). Hence we may assume that the largest number $2k \geq 6$ geodesics meet at the north pole and G contains the associated $\mathbf{Z}_{2k}$ or $D_k$. The spherical polygons have possible angles $\pi/2$, $\pi/3$, $\pi/4$, …, and at most three sides. The polygons incident to the north pole are congruent by symmetry. If they have two sides, the cone is of type (a). So we may assume that they are triangles, with sum of angles greater than $\pi$. If the nonpolar angles are both $\pi/2$, the cone is of type (b). The remaining possibilities, listing the polar (smallest) angle first, are ($\pi/3$, $\pi/3$, $\pi/2$), ($\pi/4$, $\pi/3$, $\pi/2$), and ($\pi/5$, $\pi/3$, $\pi/2$); composing symmetries yields cases (c), (d), and (e).

The symmetries of the junctions (Proposition 5.2) yield all the asserted symmetries.

*Part II: the surface looks like any tangent cone.* If a tangent cone C is a plane, then is smooth and minimal by Lemma 4.3.

Suppose that a tangent cone C is of type (a). If k = 2 and G contains a reflection, the boundary of the restriction to a small ball is weakly close to two orthogonal great circles. From reflective symmetry, the boundary must include the great circle of reflection and locally must include the totally geodesic surface of reflection. The rest of , which is minimizing and has a planar tangent cone, is smooth and minimal by Lemma 4.3. If k > 2 and G contains $D_k$, locally must include k totally geodesic surfaces of reflection, and that's all. If k > 2 and G contains $\mathbf{Z}_{2k}$, then p lies on a geodesic of fixed points. Locally, in the quotient orbifold $M_0$, one tangent cone at $p_0$ is a halfplane, and $p_0$ lies in the boundary of $_0$. Consequently locally every point of is a singular point of , and every tangent cone consists of k planes meeting at equal angles along a line tangent to . By Allard's regularity theorem 2.4, near p, consists of 2k nearly vertical minimal halfsheets meeting along the geodesic of fixed points, each with halfplane tangent cones. By Allard's boundary regularity theorem 2.4, each halfsheet is $C^1$. Consider two opposite halfsheets, $_1$ and $_2$, and any smooth flow F. Since each $_i$ is minimal, by integration by parts the first variation is just

$$\int_{_1} \dot{F} \cdot \mathbf{n}_1 + \int_{_2} \dot{F} \cdot \mathbf{n}_2 = 0,$$

because $_1 = _2$ and the conormals $\mathbf{n}_i$ are opposites. Hence each sheet is minimal as well as $C^1$ and hence $C$ (see sect. 2.4).

Suppose that a tangent cone C is of type (b). If G contains an additional reflective symmetry, locally includes the totally geodesic surface of reflection. The rest of consists of k sheets as in case (a). If G contains $D_k$, then locally includes the k totally geodesic surfaces of reflection. The rest of has a horizontal planar tangent cone and hence is smooth and minimal by Lemma 4.3. Hence we may assume that k = 2 and G contains all orientation-reserving isometries of the cube. The tangent cone is unique. Again by Allard's regularity theorem 2.4, near p, is a $C^1$ surface $C^1$ close to C. It follows that at p, is three orthogonal $C^1$ sheets, $C^1$ close to C. As in case (a) each sheet, this time decomposed into four pieces $_1$, $_2$, $_3$, $_4$, is minimal as well as $C^1$ and hence $C$ .

For the tangent cones of types (c)—(e), rich in reflective symmetry, locally simply consists of the requisite totally geodesic surfaces.

*Part III: boundary.* To classify the boundary tangent cones, note as before that by Propositions 5.1 and 5.2 the link in the unit sphere consists of geodesics meeting in a T or at equal angles in odd numbers at say the poles and in even numbers elsewhere, as described by Proposition 5.2, yielding cases (g) and (f). At other junctions, the (nontrivial) isotropy subgroups must fix or exchange the poles; the only possibility is adding the equator and horizontal reflection to case (f), yielding case (h).

Consider a point p on a smooth piece of boundary. If a tangent cone is of type (f), as in case (a) above,   locally is of type (f); smooth by Section 2.4. If a tangent cone is of type (g),   locally contains the totally geodesic halfsheet of reflection, and the rest of   is smooth and minimal by Lemma 4.3. If a tangent cone is of type (h),   locally contains the totally geodesic surface of reflection, and the rest of   is as in case (f).

*Part IV: all singularities occur.* Finally we explain that all of these cones are minimizing, say the truncated cone C  in the unit ball in $\mathbf{R}^3$. This is trivial in cases (c)—(e), where C  consists of requisite discs of reflection. We'll treat the hardest case (b); others are similar and easier. So suppose C  is of type (b) and G contains cyclic rotation $\mathbf{Z}_{2k}$ about a vertical axis and horizontal reflection. The minimizer    must contain the horizontal disc. The rest of   , in the quotient orbifold $\mathbf{R}^3$ mod $\mathbf{Z}_{2k}$, is bounded by a semicircle of longitude and the vertical axis, and hence by a projection argument must be the halfdisc, so that    = C .

Proposition 5.4 and Theorem 5.5 characterize singularities in invariant area-minimizing surfaces in 2D and 3D with oriented tangent planes (rather than choice of unit normal; see Rmk. 4.1.1). If only orientation-preserving isometries are allowed, there are no singularities by Theorem 4.1.

**5.4. Proposition (Regularity for invariant curves with oriented tangents).** *In a smooth ($C^\infty$) 2D Riemannian manifold M, among rectifiable curves oriented by unit tangent, invariant under a group G of isometries of M, with given boundary or homology, suppose that Γ minimizes length. Then on the interior Γ consists of geodesics, k of which may cross at a point with imposed reflectional (k=2) or $D_k$ (dihedral) symmetry; if k is even, each geodesic reverses orientation at the point.*

*At an isolated boundary point with multiplicity $m \geq 1$, there are the following possibilities of Figure 5.2:*

*(a) classical regularity (cf. [HaS], W1]): m geodesics meet or one geodesic changes multiplicity from $m_0 + m$ to $m_0 \geq 1$;*

*(b) a geodesic of reflective symmetry of multiplicity one and an orthogonal geodesic of multiplicity $m' \geq 2$, both of which reverse orientation when they cross; $m = 2m' - 2$;*

*(c) pairs of geodesics of multiplicities $m_i \geq 1$ meeting a geodesic $\gamma$ of reflective symmetry at angles at most $\pi/2$; the multiplicity of $\gamma$ goes from $m' \geq 0$ to 1, and $m = 2\Sigma m_i + m' - 1$.*

*All such singularities can occur.*

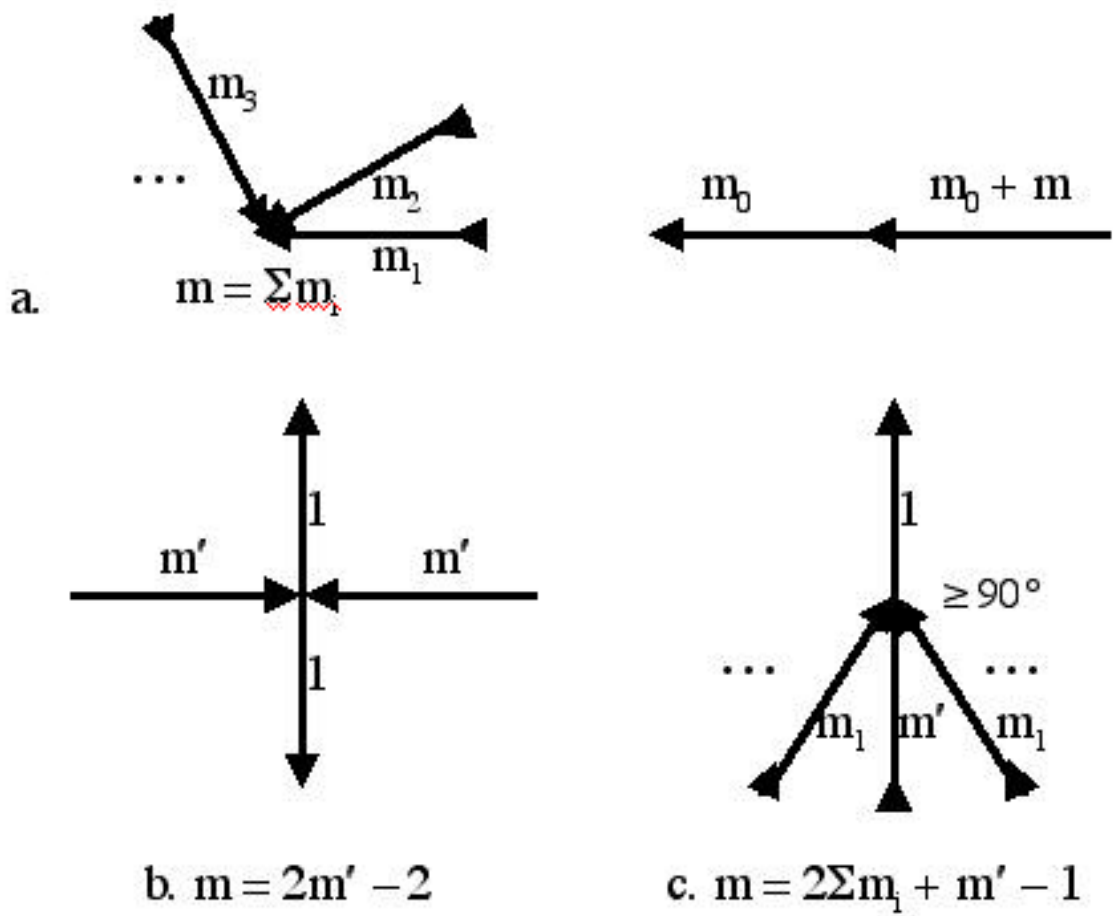

Figure 5.2
The three types of boundary singularities for curves in 2D (or surfaces in 3D), with oriented tangents, invariant under imposed (trivial or reflectional) symmetries.

*Proof.* By Proposition 5.1, it suffices to classify singular minimizing cones C in $\mathbf{R}^2$. At an interior singular point with (nontrivial) isotropy subgroup $G_1$, the tangent cone C consists of $2k \geq 4$ rays, possibly coincident, half inward and half outward. $G_1$ must be cyclic, dihedral, or generated by a single reflection. If $G_1$ is cyclic, the tangent space mod $G_1$ is a cone with angle less than $2\pi$ and C mod $G_1$ cannot be minimizing. If $G_1$ is a dihedral group [or generated by a single reflection], the tangent space mod $G_1$ is a sector with angle less than [or equal to] $\pi$ and free boundary rays of multiplicity 1/2. If there are no interior rays, the boundary rays must have opposite orientations and C consists of k geodesics with dihedral symmetry. Suppose that there are interior rays. They must have the same orientation. The boundary rays must have opposite orientation and multiplicity 1/2. Furthermore, the angle must be $\pi$. Hence C consists of two orthogonal lines with reflective symmetry.

Consider a nonclassical boundary point with (nontrivial) isotropy subgroup $G_1$. $C/G_1$ has some number of inward interior rays of multiplicity $m_i$, at least one outward boundary ray at angle at least $\pi/2$ to the interior rays, and the other boundary ray with multiplicity $m \geq \sum -1$. $G_1$ must be generated by reflection, since otherwise there can be no interior rays and hence no nonclassical examples. If $m \sum = -1$, we must have case (b); if $m \sum \geq 0$, case (c).

Consideration of the quotient orbifold shows every such singular cone in minimizing.

**5.5. Theorem (Regularity for invariant surfaces with oriented tangent planes).** *In a smooth ($C^\infty$) 3D Riemannian manifold M, among surfaces with oriented tangent planes (rectifiable currents), invariant under a group G of isometries of M, with given boundary or homology, suppose that $\Sigma$ minimizes area.*

*Then on the interior $\Sigma$ is a smoothly immersed minimal surface with only the following singularities (see Figure 5.1(a,b)):*

*(1) $k \geq 2$ "vertical" sheets of multiplicity one cross at equal angles along a smooth curve and G contains $D_k$ (dihedral) or reflection (k=2 only); if k is even, each sheet reverses orientation at the curve;*

*(2) case (1) with dihedral symmetry plus a smooth horizontal sheet; where vertical and horizontal sheets cross, both reverse orientation;*

*Along a smooth piece of boundary with multiplicity $m \geq 1$, there are the following possibilities (see Figure 5.2):*

*(a) classical regularity ([HaS], W1]): m halfsheets meet smoothly or one sheet changes multiplicity from $m_0 + m$ to $m_0 \geq 1$;*

*(b) a tangent cone consists of a plane of reflective symmetry of multiplicity one and an orthogonal plane of multiplicity $m' \geq 2$, both of which reverse orientation when they cross; $m = 2m' - 2$;*

*(c) a tangent cone consists of pairs of halfplanes of multiplicities $m_i \geq 1$ meeting a plane of reflective symmetry P at angles at most $\pi/2$; the multiplicity of P goes from $m' \geq 0$ to 1, and $m = 2\Sigma m_i + m' - 1$.*

*All such singularities can occur.*

*5.5.1. Remarks.* In boundary cases (b) and (c), a complete description remains open. The subcase of case (c) of Figure 5.3 is equivalent to a (horizontal) surface with free boundary in a (vertical) surface with boundary. As for unoriented surfaces (Remark 5.3.2), a Cantor set of singularities is possible, where the free boundary intersects the boundary. Case (b) is equivalent to a surface with both prescribed boundary and free boundary in another surface.

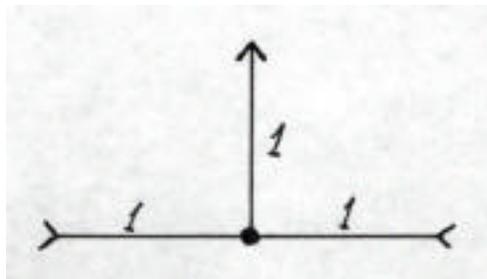

Figure 5.3
For this tangent cone of case (c) with boundary multiplicity 1, a complete description of the surface remains open. It is equivalent to free boundary in a surface with boundary.

*Proof.* The proof parallels that of Theorem 5.3. On the interior, only singularities 5.3(a,b) with $D_k$ or reflective symmetry survive; the rest cannot be oriented invariantly. Although oriented surfaces can have multiplicity, singularity 5.3(a) for example with multiplicity two (or greater) is not area minimizing, as illustrated by Figure 5.4.

Multiplicity does lead to new possibilities at the boundary, as was known already without imposed symmetry (a). Orientation stabilizes singularity 5.3(g), so that it persists beyond 90 degrees (c). Figure 5.5 shows a candidate singularity with $D_3$ symmetry which is not minimizing.

Some new arguments are required to show that the surface $\Sigma$ is diffeomorphic to its tangent cone C. Suppose that C is a plane with constant multiplicity $m \geq 1$. Possible imposed symmetries are $\mathbf{Z}_k$ (cyclic), $SO_2$, a reflection through C, and the combinations $\mathbf{D}_k$ (dihedral) and $O_2$. Without imposed symmetry, $\Sigma$ locally decomposes into m nested minimizing sheets, smooth and minimal by Lemma 4.3, coincident by the maximum principle. The same argument applies for $\mathbf{Z}_k$ and $SO_2$. If instead or in addition there is imposed reflectional symmetry, if m is even, $\Sigma$ locally decomposes into nested sheets $\Sigma_1, \ldots, \Sigma_{m/2}, \Sigma'_{m/2}, \ldots, \Sigma'_1$. Each $\Sigma_i$ is minimizing without imposed symmetry, or replacing $\Sigma_{i+} \Sigma'_i$ with $\Sigma'_{i+} \Sigma_i$ would reduce area symmetrically. Hence the m sheets are smooth coincident minimal surfaces. If m is odd, the middle sheet must lie in the totally geodesic surface of reflectional symmetry, and the previous analysis applies to the rest.

Suppose that the tangent cone C is a plane which changes multiplicity from $m_0 + m$ to $m_0 \geq 1$. Without imposed symmetry, regularity is proved by White [W1] by a decomposition. With imposed (reflectional) symmetry, consider $\Sigma$ plus a halfsheet of reflectional symmetry of multiplicity $m_0$, to decompose $\Sigma$ into $m_0$ nested smooth minimal sheets plus m simultaneously nested halfsheets with the prescribed boundary. Since all are tangent at the boundary, all coincide and lie within the surface of reflection.

Suppose that C consists of m halfplanes meeting along a line with imposed $\mathbf{Z}_k$ symmetry (so that m is a multiple of k). By the decomposition of White [W1], applied in the quotient orbifold, we may assume that m = k. Now the argument of 5.3(a) applies.

Suppose that C consists of m halfplanes meeting along a line with imposed $\mathbf{Z}_k$ symmetry. In the orbifold $M/\mathbf{Z}_k$, locally $\Sigma' = \Sigma/\mathbf{Z}_k$ is minimizing with reflectional symmetry. By the previous decomposition argument for reflectional symmetry, we may assume that m = k, and apply the argument of 5.3(a) again.

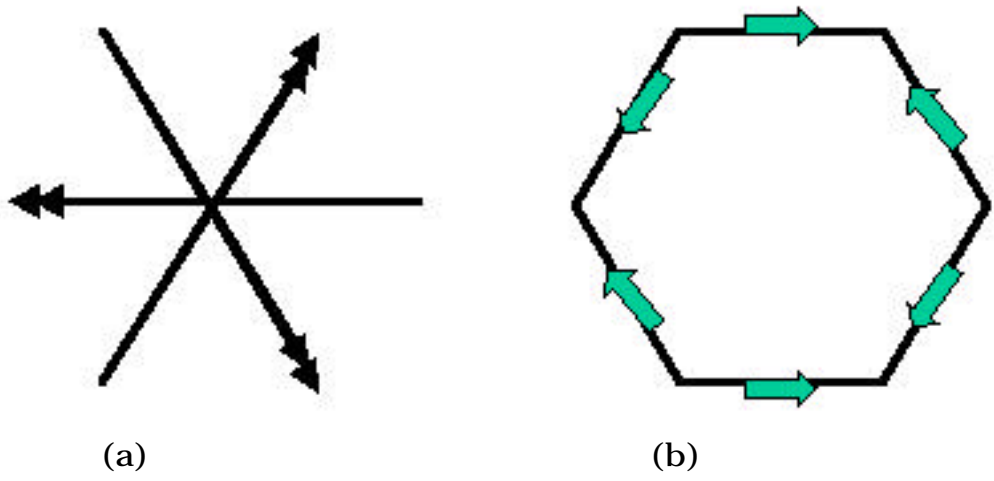

Figure 5.4
Singularity 5.3a with multiplicity two (a) is not minimizing; (b) is better.

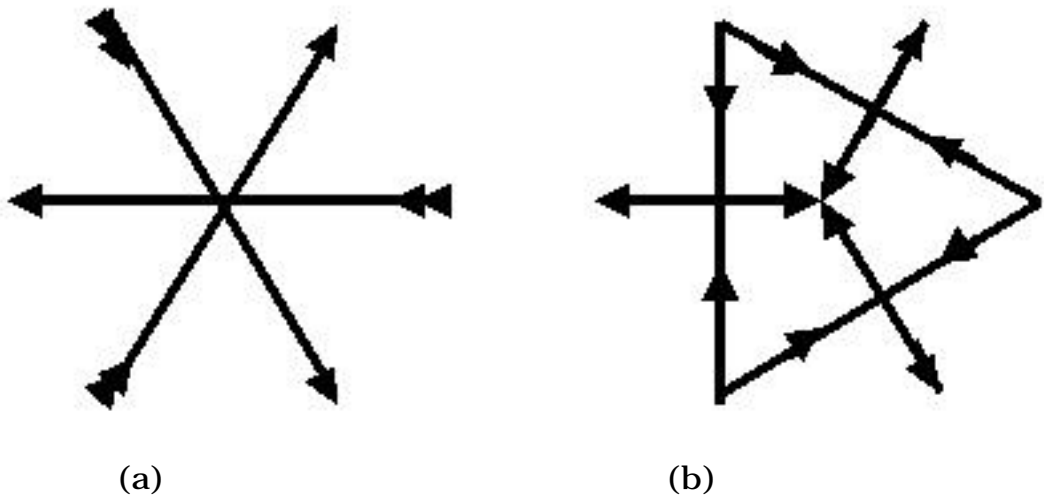

Figure 5.5
Another candidate singularity (a) which is not minimizing; (b) is better.

Proposition 5.6 and Theorem 5.7 classify singularities in invariant minimizing rectifiable curves and surfaces modulo three [M2, Chapt. 11]. The possibility of triple junctions on the interior permits new singularities (e.g. Fig. 5.8) and eliminates others (e.g. Fig. 5.7). Rectifiable curves and surfaces modulo three do not carry multiplicity, since multiplicity two is equivalent to multiplicity one with the opposite orientation. They provided early models of certain soap films [T1].

In 3D, boundary behavior remains open even for the classical case of no imposed symmetries.

**5.6. Proposition (Regularity for invariant curves modulo three).**
*In a smooth ($C^\infty$) 2D Riemannian manifold M, among rectifiable curves modulo three oriented by unit tangent, invariant under a group G of isometries of M, with given boundary or homology, suppose that Γ minimizes length.*

*At an interior point Γ may have only the following singularities of Figure 5.6:*

*(a) classical regularity: three geodesics meet at 120 degrees;*

*(b) a geodesic of reflective symmetry and an orthogonal geodesic, both of which reverse orientation when they cross;*

*(c) three geodesics cross with $D_3$ symmetry;*

*(d) k geodesics meet with $Z_k$ symmetry (k = 6, 9, ...);*

*(e) 2k geodesics meet with $D_k$ symmetry (k = 3, 6, ...);*

*At an isolated boundary point, there are the following possibilities:*

*(f) classical regularity: one geodesic, two geodesics at angle $\theta \geq 2\pi/3$, or three geodesics at $2\pi/3$;*

*(g) two geodesics at angle $2\pi/3 < \theta \leq \pi$ merge into a geodesic of reflective symmetry;*

*(h) k geodesics meet with $Z_k$ (or $D_k$) symmetry (k = 4, 5, 7, 8, ...);*

*(i) 2k geodesics meet with $D_k$ symmetry (k = 4, 5, 7, 8, ...);*

*All such singularities can occur.*

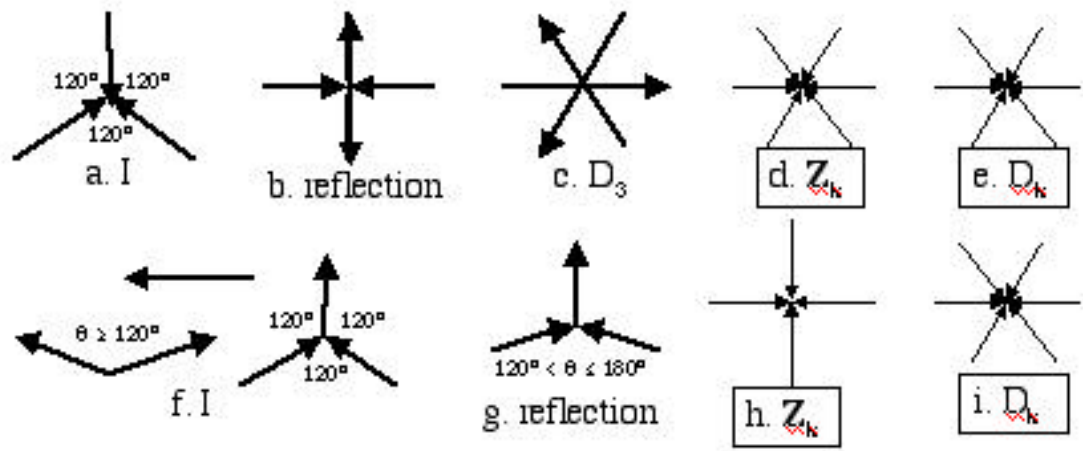

Figure 5.6
Interior (a–e) and boundary (f–i) singularities for curves in 2D (or surfaces in 3D), with oriented tangents, invariant under imposed symmetries.

*Proof.* By Proposition 5.1, it suffices to classify singular minimizing cones C in $\mathbf{R}^2$, consisting of $k_1$ outward and $k_2$ inward rays. (At an interior point, $k_2 - k_1$ is a multiple of 3.) The imposed local symmetries $G_1$ must be trivial, cyclic, dihedral, or generated by a single reflection. Without imposed symmetries, the familiar variation argument shows that two rays with the same orientation form an angle of at least $2\pi/3$, while two rays with opposite orientations must be antipodal, leading to cases (a), (f), and (g). If $G_1$ is $\mathbf{Z}_k$ (cyclic), the tangent space mod $G_1$ is a cone with angle $2\pi/k$ and C mod $G_1$ must be a single ray, yielding cases (a), (d), (f, $\theta = 180°$), and (h).

So we may assume that $G_1$ is $D_k$ (dihedral) [or generated by a single reflection]. The tangent space mod $G_1$ is a sector with angle less than [or equal to] $\pi$ and boundary rays of multiplicity 1/2.

Suppose that there are no interior rays. One boundary ray yields the previous $\mathbf{Z}_k$ cases (and $\mathbf{Z}_k \subset D_k$ suffices). Two boundary rays with the same orientation yield cases (d) and (h) and special case (f, $\theta = 180°$). Two boundary rays with opposite orientations yield cases (b) and (c); higher order examples are not minimizing, as illustrated by Figure 5.7.

Suppose that there are interior rays, along with possibly one or both boundary rays. There is room for just one interior ray, at angle at least $2\pi/3$ to a boundary ray present with the same orientation, at angle at least $\pi/2$ to a boundary ray present with opposite orientation, and at angle at least $\pi/3$ to an absent boundary ray. Since the total angle must be at least $2\pi/3$, $G_1$ is generated by a reflection. "No boundary rays" yields case (f, 120°). One boundary ray yields cases (a), (f), and (g). Two boundary rays yield case (b).

Consideration of the quotient orbifold shows every such singular cone in $\mathbf{R}^2$ minimizing.

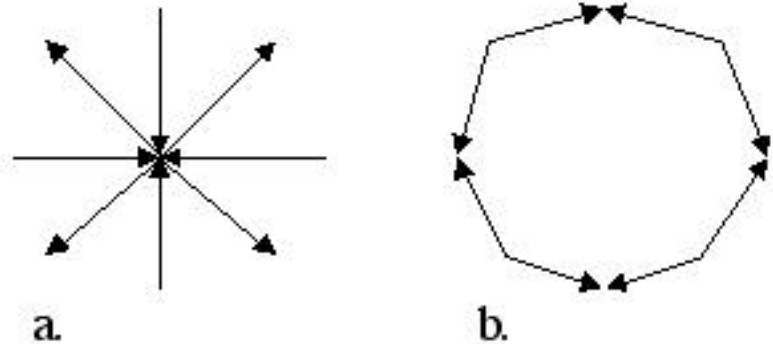

Figure 5.7
A candidate interior singularity mod 3 (a) which is not minimizing; (b) is better.

**5.7. Theorem (Regularity for invariant surfaces modulo three).**
*In a smooth ($C^\infty$) 3D Riemannian manifold M, among rectifiable currents modulo three with oriented tangent planes, invariant under a group G of isometries of M, with given boundary or homology, suppose that $\Sigma$ minimizes area.*

*Then on the interior $\Sigma$ has only the following singularities:*

*(a) three sheets meet smoothly along a smooth curve;*

*(b) two smooth minimal sheets cross along a curve with reflective symmetry;*

*(c) three smooth minimal sheets cross along a geodesic with $D_3$ symmetry;*

*(d) k smooth minimal sheets meet along a geodesic with $Z_k$ symmetry (k = 6, 9, ...);*

*(e) k smooth minimal sheets cross along a geodesic with $D_k$ symmetry (k = 3, 6, ...);*

*(f) "vertical" cases (b,c) plus a horizontal sheet, and horizontal reflection or reflection across all of the vertical sheets;*

*(g) singularity 5.3(c);*

*(h) an isolated point with tangent cone the unique minimal cone over a square pyramid with reflective symmetry as in Figure 5.8;*

*(i) an isolated point with a tangent cone over a $Z_{3k}$ or $D_{3k/2}$ symmetric link in the sphere, in which 3k geodesics from the north pole branch at triple points above the equator and symmetrically at triple points below the equator, producing 3k vertical geodesics which end at the south pole, rotated $2\pi/6k$ from those at the north pole;*

*(j) an isolated point with tangent cone with link a regular spherical tetrahedron with vertex at the north pole plus three geodesics between the poles, so that six geodesics meet at equal angles at the north pole, with imposed $D_3$ symmetry;*

*(k) an isolated point with tangent cone of Figure 5.9, with imposed $D_3$ and horizontal reflectional symmetry.*

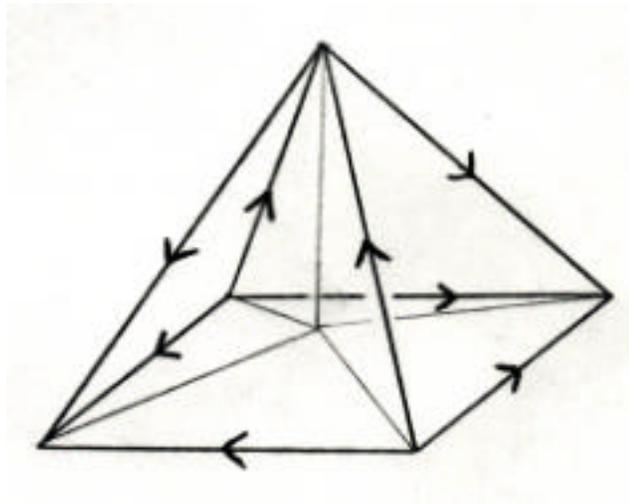

Figure 5.8
The cone over a square pyramid of case (h).

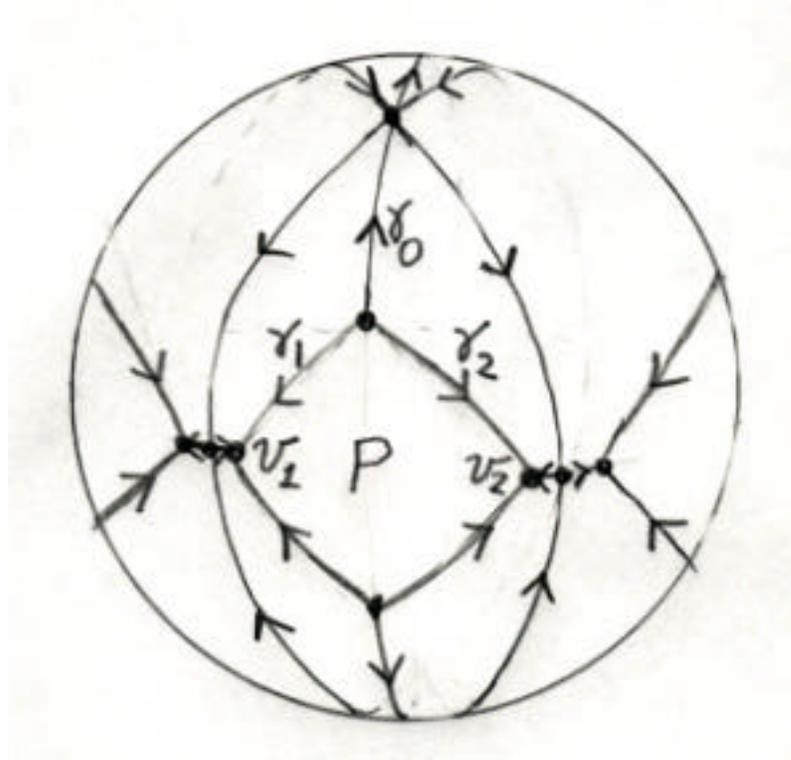

Figure 5.9
The possible tangent cone of case (k) with $D_3$ and horizontal reflectional symmetry.

*5.7.1. Remarks.* For G = I, this result is due to Taylor [T1], with higher smoothness of the singular curve by Nitsche [N]. One new singularity is the cone over the cube (case (i), k=1), which is easily seen to be minimizing if the full tetrahedral group of symmetries is imposed, and which also occurs e.g. with multiplicity three for flat chains modulo nine. I do not think that case (k) is minimizing, because each link vertex of degree four is so close to two other vertices.

We do not prove that at the isolated singularities (h,i,j) the surface is diffeomorphic to its tangent cone, although it probably follows by an argument as in Taylor [T2] or better White ([W4]; see [W2]).

Boundary behavior remains open even for the classical case of no imposed symmetries.

*Proof of Theorem 5.7. Part I: at an interior point p, a tangent cone C is a plane or of type (a)–(g).* The link the unit sphere is stationary with minimizing tangent cones. By Propositions 5.1 and 5.6, the link consists of geodesics meeting at equal angles. If there are no such junctions, C is a plane. Otherwise, we may assume that the largest number k  3 of geodesics meet at the north pole. By Gauss-Bonnet, each spherical polygon has fewer than six sides.

Suppose k = 3, so that the link consists of geodesics meeting in threes at angles of 2π/3. For orientation mod three, each spherical polygon has either two or four sides. By Gauss-Bonnet, there are either three digons, yielding case (a), or six quadrilaterals. The three quadrilaterals meeting at the north pole must be symmetric under reflection across the two sides emanating from the north pole, and hence have **Z**$_3$ symmetry. Therefore the three geodesics emanating from their lowest vertices are vertical and meet at the south pole, yielding case (g).

Suppose k = 4. If every junction has degree four, the link consists of orthogonal geodesics, and C is of type (b) or (f). Otherwise, some spherical polygon $P_1$ has $m_1 \geq 1$ vertices of degree three and $m_2 \geq 1$ vertices of degree four (one of them at the north pole), as in Figure 5.8 for example. Gauss-Bonnet says that

$$m_1(\pi/3) + m_2(\pi/2) + A = 2\pi.$$

For orientation modulo three, $m_1$ must be even. By Gauss-Bonnet, $m_1 = 2$ and $m_2 = 1$; we say that $P_1$ is of type (2,1). Call the degree-four vertex $v_0$. The adjacent polygons $P_2$, $P_3$ which share $v_0$ must likewise be of type (2,1), as must the fourth polygon $P_4$ at $v_0$. The remaining adjacent polygon $P_5$ has at least four degree-three vertices; by Gauss-Bonnet it is of type (4,0). So the cone C is a cone over a combinatorial pyramid, with at least reflective symmetry. Two opposite symmetric rays emanating from the north pole must branch at a unique height in order for the branches to meet at 120 degrees, determining the unique square pyramid, case (h).

Suppose that the singularity at the north pole is of type 5.6(c) with $D_3$ symmetry. The six polygons incident to the north pole are congruent by symmetry with $m_0 \geq 1$ angles π/3, $m_1 \geq 0$ angles π/2, and $m_2 \geq 0$ angles 2π/3. For orientation modulo three, if $m_0 = 1$, $m_2$ must be even. By Gauss-Bonnet,

$$4m_0 + 3m_1 + 2m_2 < 12,$$

and the possibilities for $(m_0, m_1, m_2)$ are (2,0,0), (2,1,0), (2,0,1), (1,2,0), (1,0,2), (1,1,0), and (1,1,2). Possibility (2,0,0) yields case (c). Possibility (2,1,0) yields case (g). For possibility (2,0,1), for orientation modulo three, the other vertex of degree six must have **Z**$_6$ symmetry, a contradiction. Possibility (1,2,0) yields case (f). For possibilities (1,0,2) and (1,1,0), the six polar triangles already have total area at least 4π, which is impossible. For the last possibility (1,1,2), some geodesic $\gamma_0$ from the north pole ends in a vertex of degree three, where it meets two other geodesics $\gamma_1$, $\gamma_2$, congruent by reflection across $\gamma_0$. The next vertex $v_1$ down $\gamma_1$ and symmetrically $v_2$ down $\gamma_2$, which are still on north polar polygons, have degree three or four.

Suppose that $v_1$ and $v_2$ have degree four as in Figure 5.10. Let $v_3$ be symmetric to $v_2$ under reflection across the next geodesic $\gamma_0$ from the pole. At $v_2$ there must be a reflectional symmetry which carries either $v_1$ or $v_3$ to the north pole, a contradiction.

Suppose that $v_1$ and $v_2$ have degree three. Then the polygon P bounded by $\gamma_1$ and $\gamma_2$ has at least three and hence four vertices of degree three. By Gauss-Bonnet there is room for at most one other vertex, of degree four, but that would break reflectional symmetry. By reflectional symmetry, P is a spherical square. Replacing it by its diagonals reduces us to the previous case $(m_1,m_2,m_3) = (1,0,2)$, and C must fall in case (k), pictured in Figure 5.9. It has more area than the area $4\varepsilon$ of that previous case, and hence compares unfavorably with a $D_3$ symmetric region on the sphere. Hence it must also have imposed horizontal reflectional symmetry.

Suppose that the singularity at the north pole is of type 5.6de with $Z_k$ or $D_{k/2}$ symmetry, $k = 6, 9, \ldots$. We may assume that there are no singularities of type 5.6(c). The k polygons incident to the north pole are congruent with $m_0 \geq 1$ angles $2\pi/k$, $m_1 \geq 0$ angles $\pi/2$, and $m_2 \geq 0$ angles $2\pi/3$ as before. For orientation modulo three, $m_0 + m_2$ must be even. By Gauss-Bonnet,

$$4m_0 + 3m_1 + 2m_2 < 12,$$

and the possibilities for $(m_0, m_1, m_2)$ are $(1,1,1)$, $(1,0,1)$, and $(1,0,3)$. For possibility $(1,1,1)$, k must be 6 (or the north polar triangles already have area at least $4\varepsilon$), yielding case (j) with $D_3$ symmetry. Possibility $(1,0,1)$ is impossible. For possibility $(1,0,3)$, the possibilities for the next tier of polygons are $(1,0,3)$, $(0,0,4)$, and $(0,1,4)$. The first yields case (i). The second yields too much area. The third violates symmetry.

The symmetries of the junctions (Proposition 5.6) yield the asserted symmetries, except for case (k), already discussed, and case (i), k=1, when the link is a spherical cube. The symmetries of the spherical cube oriented modulo three are generated by reflections across planes through antipodal edges. Reflections across any two such planes through a vertex generate the $Z_3$ about that vertex. So if G does not contain a $Z_3$, it contains reflection across at most one plane through each vertex, at most two planes total, and four faces of the cube (omitting a pair of opposite faces) provide an invariant surface of less area, a contradiction.

*Part II: the surface looks like the tangent cone.* If a tangent cone C is a plane, then  is regular by Lemma 4.3.

Suppose that a tangent cone C is of type (a). For a *new* kind of singularity, there must be some imposed symmetry, so G contains $\mathbf{Z}_3$ or G is generated by a reflection. If G contains $\mathbf{Z}_3$, the proof is like 5.3(a). So we may assume that G is generated by a reflection. In a small ball about p, replace  with a minimizer  ′ without imposed symmetry. The cheaper half of  ′, plus its reflection, is a symmetric minimizer with the same area. Hence minimizers with imposed symmetry are minimizers without imposed symmetry, and enjoy the same regularity.

Cases (b)—(g) are like cases (a)—(c) of Theorem 5.3.

In cases (h)—(k), the singularities are isolated, or there would be other points on the tangent cone with density as large as at this vertex (see [Fed1, p. 647]).

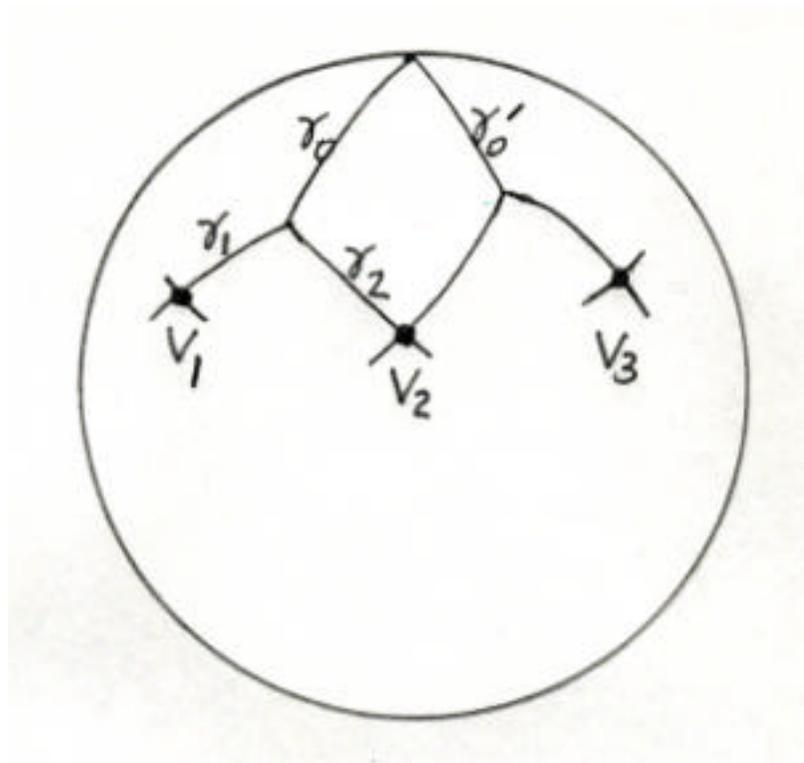

Figure 5.10
Six geodesics meeting at the north pole with $D_3$ symmetry leads to a contradiction.

**5.8. Surfaces modulo three oriented by unit normal.** If symmetries are required locally to preserve orientation of the ambient or the unit normal of the minimizing surface, then only singularities 5.7(a,d,i) survive. In this category, *all three types of singularities can occur.*

*Proof.* In $\mathbf{R}^3$, cones (a) and (d) are minimizing as products with $\mathbf{R}$ of minimizers in $\mathbf{R}^2$. Suppose that case (i) is not minimizing, and consider the minimizer bounded by a spherical cubical skeleton (or the analog for higher k). For this nice boundary, boundary regularity follows similarly to interior regularity, and near the vertices consists of three sheets meeting along a singular curve . Pass to the quotient orbifold $\mathbf{R}^3/\mathbf{Z}_3$, as in Figure 5.11. Here must separate the two halves $S_1$, $S_2$ of the sphere, and the triple curve must go from one vertex to the other. But at one vertex two of the three regions around include $S_1$, while at the other vertex two of the three regions around include $S_2$, a contradiction of the separation of $S_1$ and $S_2$.

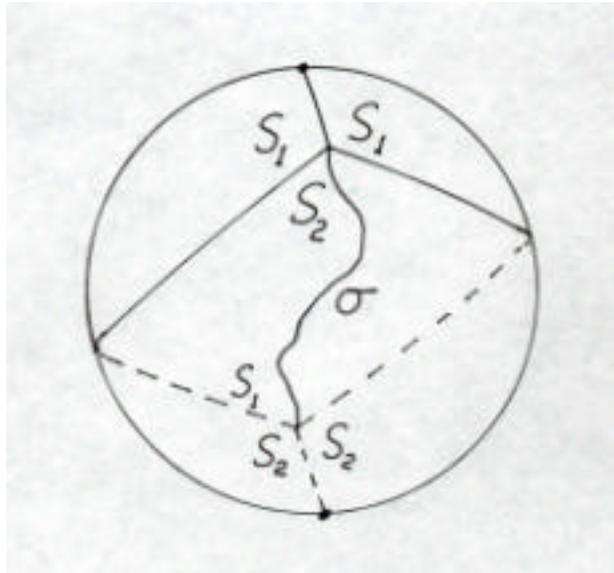

Figure 5.11
Following the purported triple curve in the quotient orbifold contradicts the fact that the surface must separate the two halves $S_1$, $S_2$ of the sphere.

**5.9. Theorem (Rectifiable currents modulo four).** *In a smooth ($C^\infty$) 3D Riemannian manifold M, among rectifiable currents modulo four oriented by unit normal, invariant under a group G of isometries of M, with given boundary or homology, suppose that Σ minimizes area. Then Σ is a smooth immersed minimal surface, which has a local mass decomposition into two area-minimizing invariant rectifiable currents modulo two. All such types of singularities occur.*

*Proof.* White [W3, 2.6] gives a local canonical mass decomposition of a flat chain modulo four into two flat chains modulo two, which is a norm-preserving bijection onto pairs $_1$, $_2$ of rectifiable currents modulo two. Because it is canonical, every symmetry in the isotropy subgroup $G_1$ either preserves or interchanges the $_i$. Let $G_2$ be the subgroup of $G_1$ of index one or two preserving the $\sum_i$. Because it is a norm-preserving bijection, $\sum$ is minimizing with imposed $G_1$ symmetry if and only if the $\sum_i$ are minimizing with imposed $G_2$ symmetry. (So actually locally S is minimizing with imposed $G_2$ symmetry.)

**5.10. Rectifiable currents modulo $v \geq 5$.** Classification of singularities for area-minimizing invariant rectifiable currents modulo 5 in 3D remains open even for $G = I$. For general orientation-preserving $G$, in 2D, geodesics can meet only $k$ at a time, with imposed $\mathbf{Z}_m$ symmetry if $k \geq 2$, where $k/m$ is an integer less than or equal to $\nu/2$.

**5.11. Soap films** (see [M2, 11.3, 13.9]). Without imposed symmetry (i.e. for $G = I$), soap films in 2D consist of geodesics meeting in threes at equal angles; and soap films in 3D consist of smooth minimal surfaces meeting in threes at equal angles along curves, which in turn meet in fours at equal angles (the cone over the tetrahedron) [T2]. For general $G$, in 2D, $k$ geodesics can meet at equal angles. Singularities in 3D include the cones over the Platonic polyhedra with full imposed symmetry, among others (as follows from consideration of the quotient orbifold).
    Even without imposed symmetry, boundary behavior remains conjectural [M2, Fig. 13.9.3].